\documentclass[reqno]{amsart}[11pt,draft,amssymb]
\usepackage{amsthm,array,amssymb,amscd,amsfonts,latexsym}

{\catcode`\@=11
\gdef\n@te#1#2{\leavevmode\vadjust{%
 {\setbox\z@\hbox to\z@{\strut#1}%
  \setbox\z@\hbox{\raise\dp\strutbox\box\z@}\ht\z@=\z@\dp\z@=\z@%
  #2\box\z@}}}
\gdef\leftnote#1{\n@te{\hss#1\quad}{}}
\gdef\rightnote#1{\n@te{\quad\kern-\leftskip#1\hss}{\moveright\hsize}}
\gdef\?{\FN@\qumark}
\gdef\qumark{\ifx\next"\DN@"##1"{\leftnote{\rm##1}}\else
 \DN@{\leftnote{\rm??}}\fi{\rm??}\next@}}
\DeclareFontFamily{OT1}{wncyr}{\hyphenchar\font45
}
\DeclareFontShape{OT1}{wncyr}{m}{n}{%
   <5> <6> <7> <8> <9> gen * wncyr
   <10> <10.95> <12> <14.4> <17.28> <20.74>  <24.88>wncyr10}{}
\DeclareFontShape{OT1}{wncyr}{m}{it}{%
   <5> <6> <7> <8> <9> gen * wncyi
   <10> <10.95> <12> <14.4> <17.28> <20.74> <24.88> wncyi10}{}
\DeclareFontShape{OT1}{wncyr}{m}{sc}{%
   <5> <6> <7> <8> <9> <10> <10.95> <12> <14.4>
   <17.28> <20.74> <24.88>wncysc10}{}
\DeclareFontShape{OT1}{wncyr}{b}{n}{%
   <5> <6> <7> <8> <9> gen * wncyb
   <10> <10.95> <12> <14.4> <17.28> <20.74> <24.88>wncyb10}{}
\input cyracc.def
\def\rus{\usefont{OT1}{wncyr}{m}{n}\cyracc\fontsize{8}{10pt}\selectfont}
\def\rusit{\usefont{OT1}{wncyr}{m}{it}\cyracc\fontsize{8}{10pt}\selectfont}

\DeclareMathSizes{9}{9}{7}{5} 

\theoremstyle{plain}

\newtheorem{theorem}{Theorem}[section]

\newtheorem{lemma}[theorem]{Lemma}

\newtheorem{corollary}[theorem]{Corollary}

\theoremstyle{definition}

\newtheorem{nothing*}[theorem]{}
\newtheorem{subnothing*}[sub]{}
\newtheorem{example}[theorem]{Example}

\theoremstyle{remark}

\def\irC{{\mathcal C}^{\rm irr}}
\def\siC{{\mathcal C}^{\rm sing}}
\def\reC{{\mathcal C}^{\rm reg}}
\def\g{{\mathfrak g}}
\def\irg{{\g}^{\rm irr}}

\def\rk{{\rm rk\,}}
\def\cod{{\rm codim}}

\begin{document}

\title[Irregular and Singular Loci
of Commuting Varieties]{Irregular and
Singular Loci of\\ Commuting Varieties}

\author[Vladimir  L. Popov]{Vladimir  L. Popov${}^*$}
\address{Steklov Mathematical Institute,
Russian Academy of Sciences, Gubkina 8,
Moscow\\ 119991, Russia}
\email{popovvl@orc.ru}

\thanks{
 ${}^*$\,Supported by Russian grants {\rus RFFI
08--01--00095a}, {\rus N{SH}--1987.2008.1},
program {\it Contemporary Problems of Theoretical
Mathematics} of the Mathe\-matics Branch of the
Russian Academy of Sciences, and by ETH,
Z\"urich.}



\dedicatory{To Bertram Kostant on the
occasion of his $80$th birthday}


\maketitle


\begin{abstract}
Let ${\mathcal C}$ be the  commuting variety of the
Lie algebra ${\mathfrak g}$ of a connected
noncommutative reductive algebraic group $G$ over an
algebraically closed field of characteristic zero.
Let $\siC$ be the singular locus of ${\mathcal C}$
and let $\irC$ be the locus  of points whose
$G$-stabilizers have dimension $>\rk\, G$. We prove
that: (a) $\siC$ is a nonempty subset of $\irC$; (b)
$\cod_{\mathcal C}\irC=5-\max\,l({\mathfrak a})$
where the maximum is taken over all simple ideals
${\mathfrak a}$ of ${\mathfrak g}$ and $l({\mathfrak
a})$ is the ``lacety'' of ${\mathfrak a}$; (c) if
${\mathfrak t}$ is a Cartan subalgebra of ${\mathfrak
g}$ and $\alpha\in {\mathfrak t}^*$ a root of
${\mathfrak g}$ with respect to ${\mathfrak t}$, then
$\overline{G({\rm Ker}\;\alpha\times {\rm
Ker}\;\alpha)}$ is an irreducible component of
$\;\irC\!$ of codimension $4$ in ${\mathcal C}$. This
yields the bound $\cod_{\mathcal C}\siC\geqslant
5-\max\,l({\mathfrak a})$ and, in particular,
$\cod{_{{}_{\mathcal C}}}\siC\geqslant 2$. The latter
 may be regarded as an evidence
in favor of the known long-standing conjecture that
${\mathcal C}$ is always normal. We also prove that
the algebraic variety ${\mathcal C}$ is rational.
\end{abstract}

\section{\bf Introduction}

\begin{nothing*} Let ${\mathfrak g}$ be a
noncommutative reductive Lie algebra
over an algebraically closed field $k$ of
characteristic zero with adjoint group $G$. Let
${\mathcal C}={\mathcal C}({\mathfrak g})$ be the
commuting variety of ${\mathfrak g}$,

\vskip -1mm

\begin{equation*}\label{com1}
{\mathcal C}={\mathcal C}({\mathfrak g}):=\{(x, y)\in
{\mathfrak g}\times {\mathfrak g}\mid [x,
y]=0\}.\end{equation*}

\vskip 2mm \noindent The known long-standing
conjectures assert that

\vskip 2mm

\begin{enumerate}
\item[(N)] ${\mathcal C}$ is normal; \item[(R)] the
ideal of regular functions on ${\mathfrak
g}\times{\mathfrak g}$ vanishing on ${\mathcal C}$ is
generated by 
functions $(a, b)\mapsto l([a, b])$ where $l\in
\mathfrak g^*$.
\end{enumerate}

\vskip 2mm

\noindent Let $\siC$ be the singular locus of
${\mathcal C}$. If Conjecture (N) is true, then,
according to the known theorem in algebraic geometry
(see, e.g., \cite[Ch.\,II, \S5, Theorem 3]{Sh}),
\begin{equation}\label{cd2}
\cod_{{\mathcal C}}\,\siC\geqslant 2.
\end{equation}
\end{nothing*}

\begin{nothing*} In this paper
we study $\cod_{\mathcal C}\,\siC$. Our approach is
based on the comparison of $\siC$ with the {\it
irregular locus} $\irC$ whose codimension we manage
to  compute. The subvariety $\irC$ is determined by
the natural action of $G$ on ${\mathcal C}$ as
follows.

Let $X$ be an algebraic variety endowed with an
action of an algebraic group $H$. For a point $x\in
X$, denote by $H(x)$ and $H_x$ respectively the
$H$-orbit and $H$-stabilizer of $x$.  If $Y$ is a
subset of $X$, then we put
\begin{gather*}
Y^{\rm reg}:=\{x\in Y\mid \dim\,H(x)\geqslant
\dim\,H(y)\hskip 2mm \mbox{\rm for every $y\in Y$}\},
\quad X^{\rm irr}:=X\setminus X^{\rm
reg}\label{reg-irr}
\end{gather*}
(although the action is not reflected in this
notation, below it is always clear from the contents
what action is meant). The set $Y^{\rm reg}$ is open
in $Y$.

As a first step we prove the following

\begin{theorem} \label{inclCC}\

\begin{enumerate}
\item[{\rm (i)}] $(0,0)\in \siC$,  so $\siC\neq
\varnothing$; \item[{\rm (ii)}] $\siC\subseteq\irC$;
\item[{\rm (iii)}] If Conjecture {\rm(R)} is true,
then $\siC=\irC$.
\end{enumerate}
\end{theorem}

\end{nothing*}

\begin{nothing*}
Then we compute $\cod_{\mathcal C}\,\irC$. To this
end we first prove that
\begin{equation}\label{bounds}
2\leqslant \cod_{\mathcal C}\,\irC\leqslant 4.
\end{equation}
We give a direct proof of \eqref{bounds} in the
framework of decomposition classes of ${\mathfrak
g}$. Actually we deduce the upper bound in
\eqref{bounds} from the following

\begin{theorem}\label{4-component}
Let ${\mathfrak t}$ be a Cartan subalgebra of
${\mathfrak g}$ and let $\alpha\in {\mathfrak t}^*$
be a root of ${\mathfrak g}$ with respect to
${\mathfrak t}$. Then
$$\overline{G({\rm Ker}\;\alpha\times {\rm
Ker}\;\alpha)},$$ is an irreducible component of
$\;\irC\!$ of codimension $4$ in ${\mathcal C}$.
\end{theorem}

Then we apply bounds \eqref{bounds} to computing
$\cod_{\mathcal C}\,\irC$. The latter problem is
immediately reduced to that for simple Lie algebras
${\mathfrak g}$. Indeed, the decomposition
\begin{equation*}\label{g-decomp}
{\mathfrak g}={\mathfrak
g}_1\oplus\cdots\oplus{\mathfrak g}_d\oplus{\mathfrak
z},
\end{equation*}
where ${\mathfrak g}_1,\ldots, {\mathfrak g}_d$ are
simple ideals and ${\mathfrak z}$ is the center of
${\mathfrak g}$, clearly, implies the decomposition
\begin{equation}\label{productcv}
{\mathcal C}({\mathfrak g})={\mathcal C}({\mathfrak
g}_1)\times\cdots\times{\mathcal C}({\mathfrak
g}_d)\times ({\mathfrak z}\times{\mathfrak z})
\end{equation}
that, in turn, implies that
\begin{equation*}\reC=
C({\mathfrak g}_1)^{\rm reg}\times\cdots\times
C({\mathfrak g}_d)^{\rm reg} \times ({\mathfrak
z}\times{\mathfrak z})
\end{equation*}
and hence
\begin{equation*}\label{mincodim}
\cod_{{\mathcal C}}\,\irC=\underset{i}{\min}\;
\cod_{{\mathcal C}({\mathfrak g}_i)}\hskip
.2mm\irC({\mathfrak g}_i).
\end{equation*}

For simple ${\mathfrak g}$, we obtain the following
complete answer:
\begin{theorem}\label{ircod}
Let ${\mathfrak g}$ be a simple Lie algebra. Then
\begin{equation*}
\cod_{\mathcal C}\,\irC=5-l,
\end{equation*}
where $l$ is the ``lacety'' of ${\mathfrak g}$, i.e.,
\begin{equation*}
l=\begin{cases}\,1 \hskip 2mm \mbox{if ${\mathfrak
g}$ is of type ${\sf A}_r$, ${\sf D}_r$, ${\sf E}_6$,
${\sf E}_7$, or
${\sf E}_8$},\\
\,2 \hskip 2mm \mbox{if ${\mathfrak g}$ is of type
${\sf B}_r$, ${\sf C}_r$, or ${\sf
F}_4$},\\
\,3 \hskip 2mm \mbox{if ${\mathfrak g}$ is of type
${\sf G}_2$}.\end{cases}
\end{equation*}
\end{theorem}

The proof of Theorem \ref{ircod}   is reduced by
means of \eqref{bounds} to finding dimensions of
certain subvarieties in the centralizers of some
nilpotent elements of some semisimple Lie algebras of
rank $\leqslant 3$. To tackle the latter problem we
go case-by-case and utilize in our arguments some
computations.
\end{nothing*}

\begin{nothing*}
The formulated results yield   the following
information about $\cod_{\mathcal C}\,\siC$. Clearly,
\eqref{productcv} implies that
\begin{equation}\label{sing}
{\mathcal C}\setminus\siC= \bigl({\mathcal
C}({\mathfrak g}_1)\setminus {\mathcal C}({\mathfrak
g}_1)^{\rm sing}\bigr)\times\cdots\times
\bigl({\mathcal C}({\mathfrak g}_d)\setminus
{\mathcal C}({\mathfrak g}_d)^{\rm sing}\bigr)\times
({\mathfrak z}\times{\mathfrak z}),\end{equation}
that, in turn, yields
\begin{equation}\label{simin}
\cod_{{\mathcal C}}\,\siC=\underset{i}{\min}\;
\cod_{{\mathcal C}({\mathfrak g}_i)}\hskip
.2mm{\mathcal C}({\mathfrak g}_i)^{\rm sing}.
\end{equation}
Thereby computing $\cod_{{\mathcal C}}\,\siC$ is
reduced to that for simple algebras ${\mathfrak g}$.
Theorems \ref{inclCC} and \ref{ircod} imply
 the following
 \begin{theorem}\label{ecs} Let ${\mathfrak g}$ be a simple
 Lie algebra. Then
 \begin{equation*}
\cod_{\mathcal C}\,\siC\geqslant 5-l,
\end{equation*}
where $l$ is the ``lacety'' of ${\mathfrak g}$.
\end{theorem}



The lower bound in \eqref{bounds} and Theorem
\ref{inclCC} show that inequality \eqref{cd2} indeed
holds for every algebra ${\mathfrak g}$. Moreover,
from Theorem \ref{ecs} and \eqref{simin} we deduce
that for some algebras ${\mathfrak g}$ a stronger
inequality holds:

\begin{corollary}\label{cor} For every
noncommutative reductive Lie algebra ${\mathfrak g}$
we have:

\begin{enumerate}
\item[{\rm(i)}]  $\cod_{\mathcal C}\hskip
.2mm\siC\geqslant 2$;

\item[\rm(ii)] If $\;{\mathfrak g}$ contains no
simple ideals of type ${\sf G}_2$, then
$\cod_{\mathcal C}\,\siC\geqslant 3$;

\item[\rm(iii)] If
$\;{\mathfrak g}$ is
simply laced, then $\cod_{\mathcal C}\,\siC\geqslant
4$.
\end{enumerate}
\end{corollary}

\begin{example} \label{example} Let ${\mathfrak g}={\mathfrak s}{\mathfrak l}_2$.
Since ${\rm SL}_2$-stabilizer of every nonzero
element of ${\mathfrak s}{\mathfrak l}_2$ is
one-dimensional, we have $\irg=\{0\}$. This implies
that $\irC=\{(0,0)\}$. By Theorem \ref{inclCC} this,
in turn, yields that $\siC=\{(0,0)\}$. As $\dim
{\mathcal C}=4$ by \eqref{dimC} below,  we have
$\cod_{\mathcal C} \siC=\cod_{\mathcal C}\irC=4$ that
agrees with Theorems \ref{ircod} and \ref{ecs}.

This case is simple enough for obtaining this and
further information directly, exploring the
equations. Namely, take an ${\mathfrak s}{\mathfrak
l}_2$-triple $e, f, h$, see \cite[\S11, 1]{Bo2}, as a
basis of ${\mathfrak g}$. Then the coordinates of
$[x_1e+x_2f+x_3h, y_1e+y_2f+y_3h]$ in this basis
generate in the polynomial algebra $k[x_1, x_2, x_3,
y_1, y_2, y_3]$ the ideal
\begin{equation}\label{sl2}
J:=(x_2y_3-x_3y_2,\; x_1y_2-x_2y_1,\;
x_1y_3-x_3y_1).\end{equation} Using \eqref{sl2} and,
e.g., a computer algebra system (we used {\sf
MAGMA}), one immediately verifies that $J$ is prime.
Hence Conjecture (R) is true in this case. This makes
it possible to prove that $\siC=\{(0,0)\}$ directly
exploring the rank of the Jacobi matrix of the system
of generators of $J$ given by~\eqref{sl2}.

Consider $x_1, x_2, x_3, y_1, y_2, y_3$ as the
standard coordinate functions on the algebra
$$
D:=\{{\rm diag}(a_1,\ldots, a_6)\in {\rm
Mat}_{6\times 6}\mid a_i\in k\}.
$$
Then ${\mathcal C}$ is a closed subset of $D$. The
group $D^*$ of invertible elements of $D$ is a
$6$-dimensional torus and \eqref{sl2} implies that
the intersection of kernels of the characters
$x_1x_2^{-1}y_1^{-1}y_2$ and $x_1
x_3^{-1}y_1^{-1}y_3$ of $D^*$ coincides with the open
subset $T:={\mathcal C}\cap D^*$ of $\;{\mathcal C}$.
Since ${\mathcal C}$ is irreducible, this means that
$T$ is a $4$-dimensional subtorus of $D^*$ and
${\mathcal C}$ is its closure in $D$. Hence
${\mathcal C}$ is the affine toric variety, namely,
the closure of $T$-orbit of the identity matrix
$I_6\in D$ with respect to the action of $\;T$ on $D$
by the left multiplication. Applying the criterion of
normality of such orbit closures
\cite[Theorem~10]{PV1}, one easily verifies that the
variety ${\mathcal C}$ is normal, i.e.,  in this case
Conjecture (N) is true as well. \quad $\square$
\end{example}

Corollary \ref{cor}(i) may be regarded as an evidence
in favor of Conjecture (N).
\end{nothing*}

\begin{nothing*}
Finally, we prove the following

\begin{theorem}\label{rational1}
The algebraic variety ${\mathcal C}$ is rational.
\end{theorem}

\end{nothing*}

\begin{nothing*}
We close this introduction by noting that the
arguments and techniques of this paper are suitable,
mutatis mutandis,  for obtaining analogous results on
commuting varieties associated with symmetric spaces
and, more generally, cyclically graded semisimple Lie
algebras, cf.\;\cite[8.5]{PV2}.
\end{nothing*}

\begin{nothing*} {\it
Notational conventions.}

\vskip 1.5mm
\begin{quote}
\begin{list}
{\renewcommand{\makelabel}{\entrylabel}} {\topsep=1mm
\leftmargin=-6mm\itemindent=0.0mm \labelwidth=1mm
\labelsep=1.3mm }
 \item{
}$x=x_s+x_n$ is the Jordan decomposition of an
element $x\in {\mathfrak g}$ with $x_s$ semisimple
and $x_n$ nilpotent.

\item{ }${\mathfrak t}$ is a Cartan subalgebra of
${\mathfrak g}$.

\item{ }$\Phi \subseteq {\mathfrak t}^*$ is the root
system of ${\mathfrak g}$ with respect to ${\mathfrak
t}$.

\item{ }$\Delta$ is a system of simple roots of
$\Phi$.

\item{ }$\Phi ^+$ is the set of positive roots of
$\Phi$ with respect to $\Delta$.

\item{ }$m:=\dim\;{\mathfrak g}$.

\item{ }$r:=\dim\;{\mathfrak t}$.

\item{ }${\mathfrak a}_\mathfrak b$ is the
centralizer of a subset $\mathfrak b$ of ${\mathfrak
g}$ in a subset $\mathfrak a$ of ${\mathfrak g}$,
\begin{equation*}\label{zentralizer}
{\mathfrak a}_{\mathfrak b}:=\{x\in \mathfrak a\mid
[x, y]=0 \text{ for all } y\in \mathfrak b\}.
\end{equation*}

\item{ }${\mathfrak z}({\mathfrak a}):={\mathfrak
a}_{\mathfrak a}$ is the center of ${\mathfrak a}$.

\item{ }${\rm Lie}\;H$ is the Lie algebra of an
algebraic group $H$.

\item{ }${\rm T}_{x}(X)$ is the tangent space to an
algebraic variety $X$ at a point $x\in X$.

\item{ }$|M|$ is the number of elements of a set $M$.

\item{ }$\langle\, S\,\rangle$ is the $k$-linear span
of a subset $S$ of a vector space over $k$.

\item{ }$I_n$ is the identity $n\times n$ matrix.

\end{list}
\end{quote}

We say that a property  holds for {\it points in
general position} of an algebraic variety $X$ if it
holds for every point laying off a proper closed
subset of $X$.

Our numeration of simple roots is that of Bourbaki
\cite{Bo1}.

\end{nothing*}

\vskip 2mm

{\it Acknowledgement.} The author thanks the referee
for remarks.

\section{\bf Morphism \boldmath$\mu$}

\begin{nothing*}

 If $x\in{\mathfrak g}$, then
 ${\rm Lie}\,G_x={\mathfrak g}_x$. Hence,
 for every $(x, y)\in {\mathfrak g}\times{\mathfrak g}$,
\begin{equation}\label{lie}
{\rm Lie}\;G_{(x, y)}={\mathfrak g}_x\cap{\mathfrak
g}_y=({\mathfrak g}_x)_y=({\mathfrak g}_y)_x.
\end{equation}

By \cite{R} we have ${\mathcal
C}=\overline{G({\mathfrak t}\times{\mathfrak t})}$.
This immediately implies that ${\mathcal C}$ is an
irreducible variety,
\begin{gather}\label{reg}
\begin{gathered}
 \reC=\{(x, y)\in {\mathcal C}\subset
{\mathfrak g}\times {\mathfrak g}\mid \dim ({\mathfrak g}_x)_y=r\}, \\
\irC= \{(x, y)\in {\mathcal C}\subset {\mathfrak
g}\times {\mathfrak g}\mid \dim ({\mathfrak
g}_x)_y>r\},
\end{gathered}
\end{gather}
and the fibers of the morphism
\begin{equation}\label{pi}
\pi_1\colon {\mathcal C}\rightarrow {\mathfrak
g},\hskip 2mm (x, y)\mapsto x,
\end{equation}
over points in general position in ${\mathfrak g}$
are isomorphic to ${\mathfrak t}$.  Since $\pi_1$ is
surjective, the latter property yields, by theorem on
dimension of fibres (see, e.g., \cite[Ch.\,I, \S6,
Theorem 7]{Sh}), that
\begin{equation}\label{dimC}
\dim\,{\mathcal C}=m+r.
\end{equation}

\end{nothing*}

\begin{nothing*}
 Consider the morphism
\begin{align}\label{mu}
\mu\colon {\mathfrak g}\times{\mathfrak
g}&\rightarrow {\mathfrak g}, \hskip 2mm (x,
y)\mapsto[x, y].
\end{align}
We have ${\mathcal C}=\mu^{-1}(0)$.

\begin{lemma} \label{lemmaker}
Let $z=(a, b)$ be a point of $\;{\mathfrak
g}\times{\mathfrak g}$. Then
\begin{enumerate}
\item[\rm (i)] $\dim\,{\rm Ker}\;
d_{z}\mu=m+\dim\,G_{z}$; \item[\rm (ii)] $\dim\,
([{\mathfrak g}, a]+[{\mathfrak g}, b])=\dim\,G(z)$.
\end{enumerate}
\end{lemma}

\vskip 2mm

\noindent {\it Proof.} We identify in the natural way
${\mathfrak g}\oplus{\mathfrak g}$ (respectively,
${\mathfrak g}$) with ${\rm T}_{z}({\mathfrak
g}\times{\mathfrak g})$ (respectively, ${\rm
T}_{\mu(z)}({\mathfrak g})$). Then \eqref{mu} implies
that the differential $d_z\mu\colon {\mathfrak
g}\oplus{\mathfrak g}\rightarrow{\mathfrak g}$ is
given by the formula
\begin{equation}\label{dm}
d_z\mu((x, y))=[x,b]+[a,y]\hskip 2mm \mbox{\rm for
every $(x, y)\in {\mathfrak g}\oplus{\mathfrak g}$}.
\end{equation}

Let $ {\mathfrak g}\times{\mathfrak g}\rightarrow k,
(x, y)\mapsto \langle x, y\rangle$, be a
 nondegenerate $G$-invariant symmetric
bilinear form  (since ${\mathfrak g}$ is reductive,
it exists). Then from \eqref{dm} we deduce that for
every $t\in {\mathfrak g}$, $(x, y)\in {\mathfrak
g}\oplus{\mathfrak g}$, we have
\begin{align*}
\langle t, d_z\mu((x, y))\rangle= \langle t, [x,
b]\rangle+\langle t, [a, y]\rangle =\langle[b, t],
x\rangle + \langle [t, a], y\rangle.
\end{align*}
This means that if ${\mathfrak g}^*$ and $({\mathfrak
g}\oplus{\mathfrak g})^*={\mathfrak
g}^*\oplus{\mathfrak g}^*$ are identified by means of
$\langle\ {,}\ \rangle$ with, respectively,
${\mathfrak g}$ and ${\mathfrak g}\oplus{\mathfrak
g}$, then the map $(d_z\mu)^*\colon {\mathfrak
g}\rightarrow {\mathfrak g}\oplus{\mathfrak g}$ dual
to $d_z\mu$
is given by the formula
\begin{equation}\label{dm*}
(d_z\mu)^*(t)= ([b, t], [t, a]).
\end{equation}
From \eqref{dm*} we deduce that
\begin{equation*}
{\rm Ker}\,(d_z\mu)^*={\mathfrak g}_a\cap{\mathfrak
g}_b,
\end{equation*}
which, together with \eqref{lie},
imply that
\begin{equation}\label{rkdzm}
\rk\,(d_z\mu)^*=m-\dim ({\mathfrak g}_a\cap{\mathfrak
g}_b) =m-\dim G_z.
\end{equation}
Now
(i) follows from \eqref{rkdzm} because
\begin{equation}\label{rkmu}
\rk\,d_z\mu=\rk\,(d_z\mu)^*. \end{equation}

Further, by \eqref{dm} the image of $d_z\mu$ is
$[{\mathfrak g},a]+[{\mathfrak g},b]$. Since $\dim
G(z)=m-\dim G_z$, this,
 \eqref{rkdzm}, and \eqref{rkmu}
imply (ii).
\quad $\square$

\vskip 2mm

\begin{lemma} \label{lemmazC}
Let $z$ be a point of $\;{\mathcal C}$. Then
\begin{enumerate}
\item[\rm(i)] $\dim\,{\rm Ker}\,d_z\mu\geqslant
\dim\,{\mathcal C}$; \item[\rm(ii)] The following
properties are equivalent:
\begin{enumerate}
\item[\rm(a)] $\dim\,{\rm Ker}\,d_z\mu=
\dim\,{\mathcal C}$;\item[\rm(b)] $z\in {\mathcal
C}^{\rm reg}$.
\end{enumerate}
\end{enumerate}
\end{lemma}
\vskip 2mm

\noindent{\it Proof.} This follows from Lemma
\ref{lemmaker}(i), \eqref{lie}, \eqref{reg}, and
\eqref{dimC}. \quad $\square$\vskip 2mm
\end{nothing*}

\begin{nothing*} {\it Proof of Theorem
{\rm \ref{inclCC}}.} By \eqref{sing}, proving (i), we
may (and shall) assume that ${\mathfrak g}$ is
simple. The point $(0,0)$ is fixed under the action
of $G$ on ${\mathfrak g}\times{\mathfrak g}$ and the
$G$-module ${\rm T}_{(0,0)}({\mathfrak
g}\times{\mathfrak g})$ is isomorphic to ${\mathfrak
g}\oplus{\mathfrak g}$. Since ${\mathfrak g}$ is
simple, this implies that every proper submodule of
${\rm T}_{(0,0)}({\mathfrak g}\times{\mathfrak g})$
is isomorphic to ${\mathfrak g}$ and, in particular,
its dimension is $m$. But ${\rm T}_{(0,0)}({\mathcal
C})$ is a nonzero submodule of ${\rm
T}_{(0,0)}({\mathfrak g}\times{\mathfrak g})$ and
$\dim\,{\rm T}_{(0,0)}({\mathcal
C})\geqslant\dim\,{\mathcal C}=m+r>m$. Hence ${\rm
T}_{(0,0)}({\mathcal C})={\rm T}_{(0,0)}({\mathfrak
g}\times{\mathfrak g})$ and therefore $\dim\,{\rm
T}_{(0,0)}({\mathcal C})=2m$.  Since $m>r$,we deduce
from here that $\dim\,{\rm T}_{(0,0)} ({\mathcal
C})>\dim\,{\mathcal C}$. This proves (i).

As ${\mathcal C}=\mu^{-1}(0)$, we have, for every
point $z\in {\mathcal C}$, the inclusion
\begin{equation}\label{TC}
{\rm T}_{z}({\mathcal C})\subseteq {\rm
Ker}\;d_z\mu.\end{equation}  Let $z\in \siC$, i.e.,
$\dim\,{\rm T}_{z} ({\mathcal C})>\dim\,{\mathcal
C}$. Then inclusion \eqref{TC} yields the inequality
$\dim {\rm Ker}\;d_z\mu>\dim\,{\mathcal C}$; whence
$z\in \irC$ by Lemma \ref{lemmazC}(ii). This proves
(ii).

Assume that Conjecture (R) is true. Then ${\rm
T}_{z}({\mathcal C})={\rm Ker}\;d_z\mu$ for every
point $z\in {\mathcal C}$ and hence the inclusion
$z\in\siC$ is equivalent to the inequality
$\dim\,{\rm Ker}\;d_z\mu>\dim\,{\mathcal C}$. By
Lem\-ma \ref{lemmazC}(ii) the latter inequality is
equivalent to the inclusion $z\in \irC$. This proves
(iii).\quad $\square$

\end{nothing*}
\begin{nothing*}{\it Remark.}
 It is claimed in \cite[Theorem 1.1]{NS}
that $\siC=\irC$ for ${\mathfrak g}={{\mathfrak
g}{\mathfrak l}}_n$. Unfortunately, the proof of this
claim given in \cite{NS} is incorrect since the
arguments on p.\,548 are based on the implicit
assumption that Conjecture~(R) is true. However,
these arguments do prove the inclusion
$\siC\subseteq\irC$ for ${\mathfrak g}={{\mathfrak
g}{\mathfrak l}}_n$ that is the particular case of
Theorem \ref{inclCC}(ii).

Similar mistake is made in paper \cite{Bre} aimed to
explore
 the singular locus of the
commuting variety of pairs  of symmetric matrices.
\end{nothing*}

\section{\bf Lower and
upper bounds for \boldmath$\cod_{\mathcal C}\,\irC$}
\begin{nothing*}
Our further arguments are based on
consideration of {\rm decomposition classes} (``{\rm
Zerle\-gungs\-klassen}'') \cite{BK}, \cite{Br1},
\cite{Br2}
(a.k.a.\,``{\rm packets}'' \cite{Sp}, \cite{E} and
``{\rm Jordan classes}'' \cite{TY}).


Recall that two elements $x, y\in {\mathfrak g}$ are
called {\it decomposition equivalent} if there is a
$g\in G$ such that ${\mathfrak g}_{x_s}={\mathfrak
g}_{g(y_s)}$ and $x_n=g(y_n)$. This defines an
equivalence relation on ${\mathfrak g}$ whose
equivalence classes are called {\it decomposition
classes}. Every decomposition class $\mathcal D$ is
an irreducible locally closed smooth $G$-stable
subvariety of ${\mathfrak g}$. All $G$-orbits in
$\mathcal D$ are of the same dimension and $\mathcal
D\subseteq \overline{\mathcal D}^{\rm reg}$. If
$z\in\mathcal D$, 
then
\begin{align}
{\mathcal
D}&=G({\mathfrak z}({\mathfrak g}_{z_s})^{\rm reg}+z_n),\notag\\
\label{dimen} \dim\,\mathcal D&=\dim {\mathfrak
z}({\mathfrak g}_{z_s})+\dim G(z).
\end{align}
\end{nothing*}

\begin{nothing*}\label{subs} As
${\mathfrak g}_{z_s}$ is a Levi subalgebra, it is
conjugate to a standard one with respect to $\Delta$.
This yields the following description of
decomposition classes.
Let $I$ be a subset of $\Delta$ (may be empty), let
$\Phi(I)$ be the set elements of $\Phi$ that are
linear combinations of elements of $I$, and let
$\Phi(I)^+:=\Phi^+\cap \Phi(I)$. Consider the Levi
subalgebra
\begin{equation*}
{\mathfrak g}(I):={\mathfrak t}\oplus \sum_{\alpha\in
\Phi(I)}{\mathfrak g}^\alpha. \end{equation*} We have
\begin{gather*}\label{zI}
{\mathfrak g}(I)={\mathfrak t}(I)\oplus{\mathfrak
s}(I) \hskip 3mm \mbox{\rm where\hskip 2mm
${\mathfrak t}(I):={\mathfrak z}({\mathfrak g}(I))$,
${\mathfrak s}(I):=[{\mathfrak g}(I), {\mathfrak
g}(I)],$}
\end{gather*}
and
\begin{gather}\label{zI}
{\mathfrak t}(I)=
\bigcap_{\alpha\in I}{\rm Ker}\;\alpha,\qquad
{\mathfrak t}(I)^{\rm reg}={\mathfrak t}(I)\setminus
\bigcup_{\alpha\notin \Phi(I)}{\rm Ker}\;\alpha.
\end{gather}

Let $x$ be a nilpotent element of ${\mathfrak s}(I)$.
 Then
\begin{equation}\label{DIxx}
\mathcal D(I, x):=G({\mathfrak t}(I)^{\rm reg} + x)
\end{equation}
is a decomposition class, every decomposition class
coincides with some $\mathcal D(I, x)$, and $\mathcal
D(I, x)=\mathcal D(J, y)$ if and only if $I=J$ and
$x$ and $y$ lay in the same orbit of the normalizer
of ${\mathfrak g}(I)$ in $G$.
In particular, since the number of nilpotent orbits
in ${\mathfrak s}(I)$ is finite, there are only
finitely many decomposition classes.

If $I$ consists of a single element $\alpha$, we
shall write ${\mathfrak t}(\alpha)$
in place of ${\mathfrak t}(\{\alpha\})$ etc.

\begin{lemma}\label{DIx} \ 
\begin{enumerate}
\item[{\rm (i)}] $\dim\mathcal D(I,
x)=m-\dim{\mathfrak s}(I)_x$; \item[{\rm (ii)}] The
following are equivalent:
\begin{enumerate} \item[\rm (a)] $\mathcal D(I, x)\subseteq
\irg$; \item[\rm (b)] $\dim{\mathfrak s}(I)_x>|I|$;
\item[\rm (c)] $x\in{\mathfrak s}(I)^{\rm irr}$.
\end{enumerate}
\end{enumerate}
\end{lemma}

\noindent{\it Proof.} If $t\in {\mathfrak t}(I)^{\rm
reg}$, then ${\mathfrak g}_t={\mathfrak g}(I)$ by
\eqref{zI}.
For $z=t+x\in \mathcal D(I, x)$, we have $t=z_s$,
$x=z_n$. This yields
\begin{equation}\label{gz}
{\mathfrak g}_{z}=({\mathfrak g}_t)_x=({\mathfrak
t}(I)\oplus{\mathfrak s}(I))_x= {\mathfrak
t}(I)\oplus{\mathfrak s}(I)_x.
\end{equation}
From \eqref{dimen} and \eqref{gz} we then deduce that
$\dim{\mathcal D}(I, x)=\dim{\mathfrak
t}(I)+m-\dim{\mathfrak g}_z=m-\dim{\mathfrak
s}(I)_x$.
This proves 
(i).

It is well known \cite{K2} that
\begin{equation}\label{irg}
\irg=\{z\in {\mathfrak g}\mid \dim {\mathfrak
g}_z>r\}.
\end{equation}
By \eqref{zI}  we have $\dim {\mathfrak t}(I)=r-|I|$,
so \eqref{gz} yields $\dim{\mathfrak
g}_z=r-|I|+\dim{\mathfrak s}(I)_x$; whence by
\eqref{irg} the
equivalence ${\rm(a)}\hskip
-1.0mm\Leftrightarrow\hskip -1.0mm{\rm(b)}$ in (ii).
Since ${\mathfrak s}(I)$ is a semisimple Lie algebra
of rank $|I|$, the equivalence ${\rm(b)}\hskip
-1.2mm\Leftrightarrow\hskip -1.2mm{\rm(c)}$ in (ii)
follows from the description of irregular loci in
reductive Lie al\-geb\-ras given by \eqref{irg}.
  \quad
$\square$

\begin{corollary}\label{cor1} Let
$\alpha\in\Delta$. Then $\mathcal D(\alpha,
x)\subseteq \irg$
if and only if $x=0$.
\end{corollary}
\noindent{\it Proof.} Since ${\mathfrak s}(\alpha)$
is isomorphic to ${\mathfrak s}{\mathfrak l}_2$ and
$({\mathfrak s}{\mathfrak l}_2)^{\rm irr}=\{0\}$, see
Example \ref{example}, the claim follows from the
equivalence ${\rm(a)}\hskip
-1.0mm\Leftrightarrow\hskip -1.0mm{\rm(c)}$ in Lemma
\ref{DIx}(ii). \quad $\square$
\end{nothing*}

 \begin{nothing*} The following statement
should be known to the experts, but I failed to find
a proper reference and shall give a short proof (for
the conjugating action of $G$ on $G$, the counterpart
of this statement is proved in \cite{St1} (see also
\cite{St2})).

\begin{lemma}\label{compirg}
\ 
\begin{enumerate}
\item[{\rm (i)}] For every $\alpha\in \Delta$, the
variety $\overline{G({\rm
Ker}\;\alpha)}=\overline{{\mathcal D}(\alpha, 0)}$ is
an irreducible component of $\irg$. \item[{\rm (ii)}]
Every irreducible component of $\;\irg$ is of this
type. \item[{\rm (iii)}] $\dim {\mathcal
D}(\alpha, 0)=m-3$ for every $\alpha\in \Delta$.
\end{enumerate}
\end{lemma}
\noindent{\it Proof.} Since in every decomposition
class all $G$-orbits are of the same dimension,
\eqref{irg} implies that $\irg$ is a union of
decomposition classes. Take a decomposition class
$\mathcal D(I, y)\subseteq {\mathfrak g}^{\rm irr}$.
Then by Lemma \eqref{DIx}(ii) we have ${\mathfrak
s}(I)\neq 0$ and $y\in {\mathfrak s}(I)^{\rm irr}$.
Hence by \cite[Theorem 5.3]{K1} there is a root
$\alpha\in I$ such that the orbit of $y$ under the
action of the adjoint group of ${\mathfrak s}(I)$
intersects the subalgebra $\sum_{\gamma\in
\Phi(I)^{+}\setminus \{\alpha\}} {\mathfrak
g}^{\gamma}$. Therefore,
\begin{equation}\label{DD}
{\mathcal D}(I, y)={\mathcal D}(I, x)\hskip
3mm\mbox{for some} \hskip 3mm
 x\in\sum_{\gamma\in \Phi(I)^{+}\setminus
\{\alpha\}} {\mathfrak g}^{\gamma}.
\end{equation}
 On
the other hand, according to \cite[5.4]{BK} (see also
\cite[39.2.2]{TY}),
\begin{equation}\label{clD}
\overline{{\mathcal D}(\alpha, 0)}\supseteq
{\mathfrak t}(\alpha)+\sum_{\gamma\in\Phi^+
\setminus\{\alpha\}}{\mathfrak g}^\gamma.
\end{equation}
From \eqref{zI}, \eqref{DIxx}, \eqref{DD},
\eqref{clD} we then deduce that
\begin{equation}\label{incl}
{\mathcal D}(I, y)\subseteq \overline{{\mathcal
D}(\alpha, 0)}=\overline{G({\rm Ker}\,\alpha)}.
\end{equation}
For every root $\gamma\in \Delta$, by Corollary
\ref{cor1} we have ${\mathcal D}(\gamma, 0)\subseteq
\irg$ and, since ${\mathfrak s}(\gamma)$ is
isomorphic to ${\mathfrak s}{\mathfrak l}_2$, Lemma
\ref{DIx}(i) yields that $ \dim {\mathcal D}(\gamma,
0)=m-3$. By virtue of \eqref{incl} this completes the
proof.\quad $\square$
\end{nothing*}

\begin{nothing*} {\it Remarks.} (a)
${\mathcal D}(\alpha, 0)={\mathcal D}(\beta, 0)$ if
and only if $W(\alpha)=W(\beta)$ where $W$ is the
Weyl group. Hence the number of irreducible
components of $\irg$ is equal to $|\Phi/W|$.

(b) Since ${\mathcal D}(\alpha, 0)= G({\mathfrak
t}(\alpha)^{\rm reg})$ and ${\mathfrak t}(\alpha)$ is
a reductive subalgebra of ${\mathfrak g}$, Lemma
\ref{compirg}(iii) is a special case of the following
more general
\begin{lemma}
Let ${\mathfrak l}$ be a reductive subalgebra of
${\mathfrak g}$ and let ${\mathfrak c}$ be its
maximal torus. Assume that ${\mathfrak c}\subseteq
{\mathfrak t}$. Then
\begin{equation*}
\dim \overline{G({\mathfrak l})}=\dim {\mathfrak c}+
|\{\alpha\in \Phi\mid
{\mathfrak c}\nsubseteq{\rm Ker}\;\alpha \}|.
\end{equation*}
\end{lemma}
\noindent{\it Proof.} Taking into account that the
image of morphism $G\times{\mathfrak l}\rightarrow
{\mathfrak g}$, $(g, x)\mapsto g(x)$, contains an
open subset of $\overline{G({\mathfrak l})}$, and the
union of maximal tori of ${\mathfrak l}$ contains an
open subset of ${\mathfrak l}$, we conclude that
points $x$ in general position in ${\mathfrak l}$ are
nonsingular points of $\overline{G({\mathfrak l})}$
and ${\rm T}_{x}(\overline{G({\mathfrak
l})})={\mathfrak l}+{\mathfrak g}(x)$. Hence $\dim
\overline{G({\mathfrak l})}=\dim({\mathfrak
l}+{\mathfrak g}(x))$.
The root decomposition of ${\mathfrak g}$ with
respect to ${\mathfrak t}$ yields
\begin{equation}\label{gx}
{\mathfrak g}(x)=\bigoplus_{\{\alpha\in \Phi\;\mid\;
{\mathfrak c}\,\nsubseteq\,{\rm Ker}\;\alpha
\}}{\mathfrak g}^\alpha.
\end{equation}
The right-hand side of \eqref{gx} is the sum of all
weight subspaces of ${\mathfrak g}$ with respect to
${\mathfrak c}$ with the nonzero weights. Since every
root space of ${\mathfrak l}$ with respect to
${\mathfrak c}$ lies in this sum, we obtain
\begin{equation*}
{\mathfrak l}+{\mathfrak g}(x)={\mathfrak
c}+\bigoplus_{\{\alpha\in \Phi\;\mid\; {\mathfrak
c}\,\nsubseteq\,{\rm Ker}\;\alpha \}}{\mathfrak
g}^\alpha;\end{equation*} whence the claim. \quad
$\square$
\end{nothing*}

\begin{nothing*}\label{notat}
It is convenient to introduce the following notation.

Let ${\mathfrak a}$ be a reductive Lie algebra and
let $x$ be its element. Put
\begin{equation}\label{ax}
{\mathcal I}({\mathfrak a}_x):= \{y\in{\mathfrak
a}_x\mid \dim ({\mathfrak a}_x)_y>{\rm
rk}\,{\mathfrak a}\}.
\end{equation}

\noindent
Clearly, ${\mathcal I}({\mathfrak a}_x)$ is empty if
and only if $x\in {\mathfrak a}^{\rm reg}$.
If $x$ is semisimple, then ${\mathfrak a}_x$ is a
reductive algebra of rank ${\rm rk}\,{\mathfrak a}$,
whence ${\mathcal I}({\mathfrak a}_x)={\mathfrak
a}_x^{\rm\hskip .7mm irr}$. For nonsemisimple $x$,
this equality, in general, does not hold (see below
Subsection \ref{G2}).
\begin{lemma}\label{calI}
Let $x\in [{\mathfrak a},{\mathfrak a}]$. Then
$\mathcal I({\mathfrak a}_x)$ is isomorphic to
${\mathfrak z}({\mathfrak a})\times\mathcal
I([{\mathfrak a},{\mathfrak a}]_x)$.
\end{lemma}
\noindent{\it Proof.} Since ${\mathfrak a}={\mathfrak
z}({\mathfrak a})\oplus[{\mathfrak a},{\mathfrak
a}]$, we have ${\mathfrak a}_x={\mathfrak
z}({\mathfrak a})\oplus [{\mathfrak a}, {\mathfrak
a}]_x$. Hence, for $z\in {\mathfrak z}({\mathfrak
a})$ and $y\in [{\mathfrak a}, {\mathfrak a}]_x$, we
have $({\mathfrak a}_x)_{z+y}={\mathfrak
z}({\mathfrak a})\oplus ([{\mathfrak a}, {\mathfrak
a}]_x)_y$. Since ${\rm rk}\,{\mathfrak a}=\dim
{\mathfrak z}({\mathfrak a})+{\rm rk}\,[{\mathfrak
a}, {\mathfrak a}]$, this and \eqref{ax} show that
$z+y\in \mathcal I({\mathfrak a}_x)$ if an only if
$y\in \mathcal I([{\mathfrak a},{\mathfrak a}])_x$;
whence the claim. \quad $\square$

\vskip 2mm

By \eqref{zI} this yields
\begin{corollary} Let $I$ be a subset of
$\Delta$ and let $x$ be an element of ${\mathfrak
s}(I)$. Then
\begin{equation}\label{dimIgx}
\dim \mathcal I({\mathfrak g}(I)_x)=r-|I|+\dim
\mathcal I({\mathfrak s}(I)_x).\end{equation}
\end{corollary}
\end{nothing*}

\begin{nothing*} It follows from \eqref{lie},
\eqref{reg}, \eqref{irg} that the restriction to
$\irC$ of the projection $\pi_1$ (see \eqref{pi})  is
a surjective morphism
\begin{equation}\label{pi1res}
\pi:=\pi_1\vert_{\irC}\colon \irC\rightarrow \irg.
\end{equation}
Let $y$ be a point of $\irg$. Then by \eqref{reg} and
\eqref{ax} we have
\begin{equation*}
\pi^{-1}(y)=\{(y, z)\in {\mathfrak g}\times{\mathfrak
g} \mid z\in\mathcal I({\mathfrak g}_y)\}.
\end{equation*}
This shows that the second projection yields an
isomorphism
\begin{equation}\label{fiberiso}
\pi^{-1}(y)\xrightarrow{\simeq} \mathcal I({\mathfrak
g}_y).
\end{equation}
\end{nothing*}

\begin{nothing*}Now we shall prove lower
bound \eqref{bounds} for $\cod_{\mathcal C}\,\irC$.

\begin{theorem}\label{>1}
$\cod_{\mathcal C}\,\irC\geqslant 2$.
\end{theorem}
\noindent{\it Proof.} Let $X$ be an irreducible
component of $\irC$. By \eqref{dimC} we have to show
that
\begin{equation}\label{cX} \dim X\leqslant
m+r-2. \end{equation} By theorem on dimension of
fibres, for every point $y\in
\pi(X)$,
we have the inequality
\begin{equation} \label{ineq}
\dim X\leqslant \dim \overline{\pi(X)} +\dim
(\pi^{-1}(y)\cap X).\end{equation}

Since $\overline{\pi(X)}$ is an irreducible variety
and ${\mathfrak g}$ is the union of decomposition
classes, there is a decomposition class $\mathcal
D=\mathcal D(I, x)$ such that
\begin{equation}\label{Dincl}
\overline{\pi(X)\cap \mathcal D}=\overline{\pi(X)}.
\end{equation}
By \eqref{Dincl} we have $\pi(X)\cap \mathcal D\neq
\varnothing$. Take a point $y\in \pi(X)\cap \mathcal
D$. Since $\pi^{-1}(y)$ is isomorphic to ${\mathcal
I}({\mathfrak g}_y)$, we have
\begin{equation}\label{pX}
\dim(\pi^{-1}(y)\cap X)\leqslant\dim {\mathcal
I}({\mathfrak g}_y).
\end{equation}

From \eqref{dimen} we obtain that
\begin{equation}\label{D=}
\dim \mathcal D= \dim {\mathfrak t}(I)+\dim G(y).
\end{equation}
It follows from \eqref{Dincl} that
$\overline{\pi(X)}\subseteq\overline {\mathcal D}$.
This and \eqref{ineq},
 \eqref{pX}, \eqref{D=} then imply that
\begin{equation}
\begin{split}
\dim X&\leqslant \dim{\mathfrak t}(I)+\dim G(y)+\dim
{\mathcal I}({\mathfrak g}_y)
\\
&=\dim{\mathfrak t}(I)+\dim G(y)+\dim {\mathfrak
g}_y-\cod_{{\mathfrak g}_y}{\mathcal I}({\mathfrak
g}_y)
\label{<=}\\
&=\dim{\mathfrak t}(I)+m-\cod_{{\mathfrak
g}_y}{\mathcal I}({\mathfrak g},
y)
\\
&\leqslant\dim{\mathfrak t}(I)+m.
\end{split}
\end{equation}

As $\pi(X)\subseteq \irg$ and all $G$-orbits in
$\mathcal D$ are of the same dimension, $\pi(X)\cap
\mathcal D\neq \varnothing$ implies that $\mathcal
D\subseteq \irg.$
 By Lemma \ref{DIx} this yields
 $I\neq\varnothing$, hence
  $\dim {\mathfrak t}(I)\leqslant r-1$.

If $\dim {\mathfrak t}(I)\leqslant
  r-2$, then \eqref{<=} implies
  \eqref{cX}.

  So it remains to consider
  the case where $\dim {\mathfrak t}(I)=r-1$, i.e.,
  $\mathcal D=\mathcal D(\alpha, x)$
  for some root $\alpha$.
  Since $\mathcal D\subseteq
  \irg$,  we deduce from
  Corollary \ref{cor1}
  that $x=0$. By \eqref{DIxx} this means that
  $y$ is a semisimple element.
  Hence ${\mathcal I}({\mathfrak g}_y)=({\mathfrak g}_y)^{\rm
  irr}$. By Lemma~\ref{compirg} this
  implies that $\cod_{{\mathfrak g}_y}
  {\mathcal I}({\mathfrak g}_y)=3$. Plugging this in \eqref{<=},
  we obtain $\dim X\leqslant m+r-4$;
  whence \eqref{cX}. This completes the
  proof. \quad $\square$
  \end{nothing*}
\begin{nothing*} The following statement will be
used in the proof of upper bound \eqref{bounds} for
$\cod_{\mathcal C}\,\irC$.
\begin{lemma} \label{fibers} Let
$\varphi\colon X\rightarrow Y$ be a dominant morphism
of algebraic varieties. Assume that $Y$ is
irreducible {\rm(}but $X$ may be not\hskip
.2mm{\rm)}. Then
\begin{enumerate}
\item[\rm(i)] There are an integer $c\geqslant 0$ and
an irreducible component $Z$ of $X$ such that
\begin{enumerate}
\item[\rm (a)] $\dim \varphi^{-1}(y)=c$ for points
$y$ in general position in $Y$; \item[\rm (b)] $\dim
Z=c+\dim Y$; \item[\rm (c)]
$\overline{\varphi(Z)}=Y$.
\end{enumerate}
\item[\rm(ii)] If the fibers of $\varphi$ over points
in general position in $Y$ are irreducible, then $Z$
is the unique irreducible component of $Y$ whose
image under $\varphi$ is dense in $Y$ and there is an
open subset $U$ of $Y$ such that $Z=\overline
{\varphi^{-1}(U)}$.
\end{enumerate}
\end{lemma}
\noindent{\it Proof.}
Since  $\varphi$ is dominant and  $Y$ is irreducible,
there is an irreducible component of $X$ whose image
under $\varphi$ is dense in $Y$. Let $Z_1, \ldots,
Z_n$ be all such components. Put
$\psi_i:=\varphi\vert_{Z_i}\colon
Z_i\rightarrow\nobreak Y$. By theorem of dimension of
fibers applied to $\psi_i$ there is an integer
$c_i\geqslant 0$ such that
 $\dim\psi_i^{-1}(y)=c_i$
 for points $y$ in general
 position in $Y$ and
\begin{equation}\label{Zdc}
\dim Z_i=c_i+\dim Y.\end{equation}

 Put $c=\underset{i}{\max}\, c_i$ and let
$c=c_{i_0}$. By construction, for points $y$ in
general position in $Y$, we~have
\begin{equation}\label{phi-1y}
\varphi^{-1}(y)=\bigcup_{i=1}^n \psi_i^{-1}(y).
\end{equation}
Hence $\dim \varphi^{-1}(y)=c$. This and \eqref{Zdc}
show that we can take $Z:=Z_{i_0}$. This proves (i).

If  $\varphi^{-1}(y)$ in \eqref{phi-1y}
 is irreducible, then
 $\varphi^{-1}(y)=\psi^{-1}_{i_0}(y)$
 because of the dimension reason.
 Hence there is an open subset $U$ in $Y$ such
 that for every $i$ we have
 $\psi^{-1}_i(U)\subseteq Z$. Since $Z_i$
 is irreducible, we obtain
 $Z_i=\overline{\psi^{-1}_i(U)}\subseteq
 Z$, i.e., $i=i_0$. This proves
 (ii).
 \quad $\square$
\end{nothing*}
\begin{nothing*} {\it Proof of Theorem
{\rm \ref{4-component}}.} Let $Y$ be an irreducible
component of $\irg$, let $X:=\pi^{-1}(Y)$ (see
\eqref{pi1res}), and let $\varphi:=\pi\vert_{X}\colon
X\rightarrow Y$. By Lemma \ref{fibers} there is an
irreducible component $Z$ of $X$ and an integer
$c\geqslant 0$ such that properties (a), (b), (c) in
the formulation of this lemma hold.

The variety $Z$ is an irreducible component of
$\irC$. Indeed, since $Z$ is irreducible, there is an
irreducible component $Z'$ of $\irC$ containing $Z$.
By
(c) we have $Y\subseteq\overline{\pi(Z')}$. Since $Y$
is an irreducible component of $\irg$ and
$\overline{\pi(Z')}$ is an irreducible subvariety of
$\irg$, this implies $Y=\overline{\pi(Z')}$. Hence
$Z'\subseteq X$. Since $Z$ is a maximal irreducible
closed subset of $X$ and $Z\subseteq Z'$, we have
$Z=Z'$.

By Lemma \ref{compirg} there is $\alpha\in \Delta$
such that
\begin{equation}\label{Y}
Y=\overline{G({\rm Ker}\;\alpha)}.
\end{equation}
For every $y\in({\rm Ker}\;\alpha)^{\rm reg}$, we
have ${\mathfrak g}_y={\mathfrak g}(\alpha)$. Since
${\mathfrak s}(\alpha)$ is isomorphic to ${\mathfrak
s}{\mathfrak l}_2$, this yields ${\mathfrak g}_y^{\rm
irr}={\rm Ker}\;\alpha$. As $y$ is semisimple, the
discussion in Subsection {\rm \ref{notat}} then
implies that
\begin{equation}\label{fiber}
\varphi^{-1}(y)=\{(y, z)\in {\mathcal C}\mid z\in
{\rm Ker}\;\alpha\}.
\end{equation}
This, in particular, shows that $\varphi^{-1}(y)$ is
irreducible and
\begin{equation}\label{r-1}
c=\dim \varphi^{-1}(y)=r-1. \end{equation}

From \eqref{Y}, \eqref{fiber}, and Lemma \ref{fibers}
it clearly follows that $Z=\overline{G({\rm
Ker}\;\alpha\times {\rm Ker}\;\alpha)}$. Since by
Lemma \ref{compirg} we have
\begin{equation*}\label{m-3}
\dim Y=m-3,
\end{equation*} Lemma \ref{fibers}
and \eqref{r-1} yield that $\dim Z=m+r-4$. By virtue
of \eqref{dimC} this completes the proof. \quad
$\square$
\end{nothing*}
\begin{nothing*}
Theorems \ref{4-component} and \ref{>1} reduce
computing $\cod_{\mathcal C}\irC$ to finding out
whether there are irreducible components of $\irC$ of
codimensions 2 and 3 or not. We now turn to solving
this problem.
\begin{lemma} \label{predecomp}
Let $C$ be an irreducible component of $\irC$. Then
there is a decomposition class $\mathcal D\subseteq
\irg$ such that $C$ is the closure of one of the
irreducible components of $\pi^{-1}(\mathcal D)$.
\end{lemma}
\noindent{\it Proof.} Since $\irg$ is the union of
decomposition classes and $\pi$ is surjective, we
have
\begin{equation*}\label{unioncomp}
\irC=\bigcup_{i=1}^{n}C_i,
\end{equation*}
where $C_1,\ldots, C_n$ is the set of all irreducible
components of all varieties $\pi^{-1}(\mathcal D)$
where $\mathcal D$ runs through the set of all
decomposition classes contained in $\irg$. Hence
\begin{equation*}
 C=\bigcup_{i=1}^{n}\overline{C\cap C_i}.
\end{equation*}
Since $C$ is irreducible, this implies that
$C=\overline{C\cap C_j}$  for some $j$. Hence
$C\subseteq \overline{C_j}\subseteq \irC$. As
$\overline{C_j}$ is irreducible and $C$ is a maximal
irreducible closed subset of $\irC$, from this we
deduce that $C=\overline{C_j}$. \quad $\square$

\begin{lemma}\label{preimD}\
\begin{enumerate}
\item[\rm(i)] There is an irreducible component $Z$
of $\pi^{-1}(\mathcal D(I, x))$ such that
\begin{equation}\label{codimCZ}
\cod_{\mathcal C} Z=\cod_{{\mathfrak s}(I)_x}\mathcal
I({\mathfrak s}(I)_x) + |I|
\end{equation}
and $\pi(Z)$ is dense in $\mathcal D(I, x)$.
\item[\rm(ii)] $\cod_{\mathcal C} Z'\geqslant
\cod_{\mathcal C} Z$ for every irreducible component
$Z'$ of $\pi^{-1}(\mathcal D(I, x))$.
\end{enumerate}
\end{lemma}
\noindent{\it Proof.} Let $z$ be a point of
${\mathfrak t}(I)^{\rm reg}$ and let $y=z+x$. We have
${\mathfrak g}_z={\mathfrak g}(I)$; whence
${\mathfrak g}_y=({\mathfrak g}_z)_x={\mathfrak
g}(I)_x$. From this and \eqref{fiberiso} we deduce
that $\pi^{-1}(y)$ is isomorphic to $\mathcal
I\bigl({\mathfrak g}(I)_x\bigr)$. Therefore, by
\eqref{dimIgx} we have
\begin{equation}\label{mid}
\dim \pi^{-1}(y)=r-|I|+\dim \mathcal I({\mathfrak
s}(I)_x).
\end{equation}

It follows from \eqref{DIxx} and \eqref{mid} that
dimension of fiber of $\pi$ over every point of
$\mathcal D(I, x)$ is equal to $r-|I|+\dim \mathcal
I({\mathfrak s}(I)_x)$. By theorem on dimension of
fibers this, Lemma \ref{DIx}(i), and Lemma
\ref{fibers} imply that there is an irreducible
component $Z$ of $\pi^{-1}(\mathcal D(I, x))$ such
that $\pi(Z)$ is dense in $\mathcal D(I, x)$,
\begin{equation}\label{dimbase}
\dim Z=m-\dim{\mathfrak s}(I)_x+r-|I|+\dim\mathcal
I({\mathfrak s}(I)_x),
\end{equation}
and $\dim Z\geqslant \dim Z'$ for  every irreducible
component $Z'$ of $\pi^{-1}(\mathcal D(I, x))$. Since
\eqref{codimCZ} follows from \eqref{dimbase} and
\eqref{dimC}, this completes the proof. \quad
$\square$
\end{nothing*}

\begin{nothing*}\label{explain}
Theorems \ref{4-component}, \ref{>1} and Lemmas
\ref{predecomp}, \ref{preimD}, \ref{DIx}(ii) reduce
our problem to finding the numbers
\begin{equation}\label{cIx} c(I,
x):=\cod_{{\mathfrak s}(I)_x}\mathcal I({\mathfrak
s}(I)_x)+|I|
\end{equation}
for all the cases where
\begin{equation}\label{nilpir}
1\leqslant |I|\leqslant 3 \hskip 2.5mm\mbox{and $x$
is a nilpotent element of ${\mathfrak s}(I)^{\rm
irr}$}.
\end{equation}
Namely, $\cod_{\mathcal C}\irC=2$ if and only if
there is a subset $I$ of $\Delta$ and a nilpotent
element $x\in {\mathfrak s}(I)^{\rm irr}$ such that
$|I|\leqslant 2$ and $c(I, x)=2$. If $\cod_{\mathcal
C}\irC\neq 2$, then $\cod_{\mathcal C}\irC=3$ if and
only if there is a subset $I$ of $\Delta$ and a
nilpotent element $x\in {\mathfrak s}(I)^{\rm irr}$
such that $|I|\leqslant 3$ and $c(I, x)=3$. If
$\cod_{\mathcal C}\irC\neq 2$ and $3$, then
$\cod_{\mathcal C}\irC=4$ by Theorems \ref{>1} and
\ref{4-component}.
\end{nothing*}

\begin{nothing*}
If $x=0$, then ${\mathfrak s}(I)_x={\mathfrak s}(I)$
and $\mathcal I({\mathfrak s}(I)_x)={\mathfrak
s}(I)^{\rm irr}$. So by Lemma \ref{compirg} we have
\begin{equation}\label{x=0}
c(I, 0)=3+|I|.
\end{equation}
\end{nothing*}

\begin{nothing*}\label{I1} This covers the
cases where $|I|$=1. Indeed, then ${\mathfrak s}(I)$
is isomorphic to ${\mathfrak s}{\mathfrak l}_2$,
hence $x=0$ by \eqref{nilpir} and therefore by
\eqref{x=0} in this case we have
$$c(I, x)=4.$$
\end{nothing*}

\begin{nothing*}
To explore the cases $|I|=2$ and $3$, in the next
section we obtain  a necessary information on
$\cod_{{\mathfrak a}_x}\mathcal I({\mathfrak a}_x)$
for some semisimple Lie algebras ${\mathfrak a}$ of
rank $\leqslant 3$ and nonzero nilpotent elements
$x\in {\mathfrak a}^{\rm irr}$.

Below, for such ${\mathfrak a}$, we denote  by $\Psi$
the root system of ${\mathfrak a}$ with respect to a
fixed Cartan subalgebra ${\mathfrak c}$ and by
$\{\alpha_i\}$ a system of simple roots of $\Psi$. We
fix
 a
Chevalley system $(X_{\alpha})_{\alpha\in \Psi}$ of
$({\mathfrak a}, {\mathfrak c})$ and put
$H_\alpha=[X_{-\alpha}, X_\alpha]$, cf.\;\cite[\S2,
4]{Bo2}. For classical ${\mathfrak a}$, we take
$X_\alpha$ and $H_\alpha$ as in \cite[\S13]{Bo2}. The
integers $N_{\alpha, \beta}$  for $\alpha, \beta,
\alpha+\beta\in \Psi$ are defined by the equality
$[X_\alpha, X_\beta]=N_{\alpha,
\beta}X_{\alpha+\beta}$, cf.\;\cite[\S2, 4]{Bo2}.
\end{nothing*}

\section{\bf \boldmath$\mathcal I({\mathfrak a}_x)$
for some algebras ${\mathfrak a}$ of rank
\boldmath$\leqslant 3$}

\begin{nothing*} \label{sl2sl2} {\it Case
${\mathfrak a}={\mathfrak s}{\mathfrak
l}_2\oplus{\mathfrak s}{\mathfrak l}_2$.}

Up to an outer automorphism, $x=(y, 0)$ for a nonzero
nilpotent element $y\in{\mathfrak s}{\mathfrak l}_2$.
Then ${\mathfrak a}_x= \langle y\rangle\oplus
{\mathfrak s}{\mathfrak l}_2$ and so $\mathcal
I({\mathfrak a}_x)=\langle y\rangle\oplus\{0\}$.
Therefore,
\begin{equation*}\label{sl2sl2} \cod_{{\mathfrak a}_x}\mathcal
I({\mathfrak a}_x)= 3.
\end{equation*}
\end{nothing*}

\begin{nothing*}\label{sl3} {\it Case
${\mathfrak a}={\mathfrak s}{\mathfrak l}_3$.}

In this case, the subregular orbit is the unique
nonzero nilpotent orbit of the adjoint group of
${\mathfrak a}$ in ${\mathfrak a}^{\rm irr}$. It
contains $X_{\alpha_1}\!$. Since ${\rm Ker}\,\alpha_1
=\langle H_{\alpha_1}+2H_{\alpha_2}\rangle$ and, for
every $\alpha\in \Psi$, the subalgebra ${\mathfrak
a}_{X_{\alpha}}\!$ is the linear span of ${\rm
Ker}\,\alpha$
 and all the
$X_{\beta}$'s such that $\alpha+\beta\notin \Psi$,
this yields
\begin{equation}\label{centsl3}
{\mathfrak a}_{X_{\alpha_1}}\hskip -1mm= \langle\,
X_{\alpha_1}, X_{\alpha_1+\alpha_2}, X_{-\alpha_2},
H_{\alpha_1}+2H_{\alpha_2}\rangle.
\end{equation}

As $N_{\alpha_1+\alpha_2, -\alpha_2}=-1$, we obtain,
for $a, b, c, d\in k$ and $y=aX_{\alpha_1}+
bX_{\alpha_1+\alpha_2}+ cX_{-\alpha_2}+
d(H_{\alpha_1}+2H_{\alpha_2})\in {\mathfrak
a}_{x_{\alpha_1}}\!$, that
\begin{equation}\label{expsl3}
\begin{aligned}
{[}y, X_{\alpha_1+\alpha_2}{]}
&=cX_{\alpha_1}\!+3dX_{\alpha_1+\alpha_2},\\
[y, X_{-\alpha_2}]&=-bX_{\alpha_1}\!-3d X_{-\alpha_2},\\
[y, H_{\alpha_1}+2H_{\alpha_2}]&=
3bX_{\alpha_1+\alpha_2}-3cX_{-\alpha_2}.
\end{aligned}
\end{equation}
From \eqref{centsl3} and \eqref{expsl3} we deduce
that
\begin{equation*}\label{dimcentrsl3}
\begin{aligned}
\dim ({\mathfrak a}_{X_{\alpha_1}})_{y}& =\dim
{\mathfrak a}_{X_{\alpha_1}}-\dim[y, {\mathfrak
a}_{X_{\alpha_1}}]
\\
&=4-{\rm rk} A,
\end{aligned}
\end{equation*}
where
$$A=\begin{bmatrix}c&3d&0\\
-b&0&-3d\\
0&3b&-3c\end{bmatrix}.$$

Since ${\rm rk}\,A\leqslant 2$ for all $b, c, d$, and
${\rm rk}\,A< 2$ only for $b=c=d=0$, this implies
that $\mathcal I({\mathfrak a}_{X_{\alpha_1}})$ is
the center of ${\mathfrak a}_{X_{\alpha_1}}\!$,
\begin{equation}\label{centersl3}
\mathcal I({\mathfrak a}_{X_{\alpha_1}})=
\langle\,X_{\alpha_1}\,\rangle.
\end{equation}
Therefore, by \eqref{centsl3} and \eqref{centersl3}
we have
$$
\cod_{{\mathfrak a}_{X_{\alpha_1}}}\mathcal
I({\mathfrak a}_{X_{\alpha_1}})= 3.
$$
\end{nothing*}

\begin{nothing*}\label{so51} {\it Case
${\mathfrak a}={\mathfrak s}{\mathfrak o}_5$.}

Like in the previous case we obtain
\begin{equation}\label{cso5}
{\mathfrak a}_{X_{\alpha_2}}=\langle\, X_{\alpha_2},
X_{\alpha_1+2\alpha_2}, X_{-\alpha_1},
H_{\alpha_1+\alpha_2}\,\rangle.
\end{equation}
By the dimension reason \eqref{cso5} implies that
$X_{\alpha_2}$ is a subregular element of ${\mathfrak
a}$. If $a, b, c, d\in k$, then, for the element
$y=aX_{\alpha_2}\! +bX_{\alpha_1+2\alpha_2}\!+
cX_{-\alpha_1}\!+dH_{\alpha_1+\alpha_2}\in {\mathfrak
a}_{X_{\alpha_2}}$, we have
\begin{equation}\label{zso5}
\begin{aligned}
{[}y, X_{\alpha_1+2\alpha_2}{]}&=
2dX_{\alpha_1+2\alpha_2},\\
{[}y, X_{-\alpha_1}{]}&=
-2dX_{-\alpha_1},\\
[y, H_{\alpha_1+\alpha_2}]&=
2bX_{\alpha_1+2\alpha_2}- 2cX_{-\alpha_1}.
\end{aligned}
\end{equation}
From \eqref{cso5} and \eqref{zso5} we deduce that
\begin{equation*}
\begin{aligned}
\dim ({\mathfrak a}_{X_{\alpha_2}})_{y}& =\dim
{\mathfrak a}_{X_{\alpha_2}}-\dim[y, {\mathfrak
a}_{X_{\alpha_2}}]
\\
&=4-{\rm rk}\,A,
\end{aligned}
\end{equation*}
where
$$A=\begin{bmatrix}2d&0\\
0&-2d\\
2b&-2c\end{bmatrix}.$$

 Since ${\rm rk}\,A=2$
if $d\neq 0$, and ${\rm rk}\,A< 2$ otherwise, this
implies that
\begin{equation}\label{so5centersl3}
\mathcal I\bigl({\mathfrak a}_{X_{\alpha_2}}\bigr)=
\langle\,X_{\alpha_1}, X_{2\alpha_1+\alpha_2},
X_{-\alpha_2}\,\rangle.
\end{equation}
Therefore, by \eqref{cso5} and \eqref{so5centersl3}
we have
$$
\cod_{{\mathfrak a}_{X_{\alpha_2}}}\mathcal
I({\mathfrak a}_{X_{\alpha_2}})= 1.
$$
\end{nothing*}

\begin{nothing*}\label{so52}
There is a unique nonzero nilpotent orbit $\mathcal
O$ of the adjoint group of ${\mathfrak a}$ distinct
from the subregular one, see, e.g., \cite{CM}. Its
dimension is $4$. Since
\begin{equation}\label{nonsubreg}
{\mathfrak a}_{X_{\alpha_1}}= \langle X_{\alpha_1},
X_{\alpha_1+\alpha_2}, X_{\alpha_1+2\alpha_2},
X_{-\alpha_2}, X_{-\alpha_1-2\alpha_2},
H_{\alpha_1+2\alpha_2} \rangle,
\end{equation}
by the dimension reason we have
$X_{\alpha_1}\in\mathcal O$. As
$N_{\alpha_1+\alpha_2, -\alpha_2}=2$,
$N_{\alpha_1+2\alpha_2, -\alpha_2}=-1$,
$N_{\alpha_1+\alpha_2, -\alpha_1-2\alpha_2}=1$, we
obtain, for $a, b, c, d, e, f\in k$ and $y=
aX_{\alpha_1}+ bX_{\alpha_1+\alpha_2}+
cX_{\alpha_1+2\alpha_2}+ dX_{-\alpha_2}+
eX_{-\alpha_1-2\alpha_2}+ fH_{\alpha_1+2\alpha_2}\in
{\mathfrak s}(I)_{X_{\alpha_1}}$, that
\begin{equation}\label{zso55}
\begin{aligned}
{[}y, X_{\alpha_1+\alpha_2}{]}&=
-2dX_{\alpha_1}-eX_{-\alpha_2}+
fX_{\alpha_1+\alpha_2},\\
[y, X_{\alpha_1+2\alpha_2}]&= dX_{\alpha_1+\alpha_2}+
eH_{\alpha_1+2\alpha_2}+
2fX_{\alpha_1+2\alpha_2},\\
[y, X_{-\alpha_2}]&= 2bX_{\alpha_1}-cX_{\alpha_1+
\alpha_2}-fX_{-\alpha_2},\\
[y, X_{-\alpha_1-2\alpha_2}]&
=bX_{-\alpha_2}-cH_{\alpha_1+2\alpha_2}
-2fX_{-\alpha_1-2\alpha_2},\\
[y, H_{\alpha_2+2\alpha_1}]&= -bX_{\alpha_1+\alpha_2}
-2cX_{\alpha_1+2\alpha_2}
+dX_{-\alpha_2}+2eX_{-\alpha_1-2\alpha_2}.
\end{aligned}
\end{equation}
It follows from \eqref{nonsubreg} and \eqref{zso55}
that
\begin{equation*}
\begin{aligned}
\dim ({\mathfrak a}_{X_{\alpha_1}})_{y}& =\dim
{\mathfrak a}_{X_{\alpha_1}}\!\!-\dim[y, {\mathfrak
a}_{X_{\alpha_1}}]
\\
&=6-{\rm rk}\,A,
\end{aligned}
\end{equation*}
where
$$A=\begin{bmatrix}-2d&f&0&-e&0&0\\
0&d&0&0&2f&e\\
2b&-c&0&-f&0&0\\
0&0&0&b&-2f&0\\
0&b&-2c&d&2e&0
\end{bmatrix}.$$

An elementary exploration of ${\rm rk}\,A$ as
function of $b, c, d, e, f$ yields that $\mathcal
I({\mathfrak a}_{x_{\alpha_1}})$ is the union of six
$3$-dimensional linear subspaces of ${\mathfrak
a}_{X_{\alpha_1}}$:
\begin{equation}\label{333so5}
\begin{aligned}
\mathcal I({\mathfrak a}_{X_{\alpha_1}})&= \langle
X_{\alpha_1}, X_{\alpha_1+2\alpha_2},
X_{-\alpha_2}\rangle \cup \langle X_{\alpha_1},
X_{\alpha_1+\alpha_2},
X_{\alpha_1+2\alpha_2}\rangle\\
&\hskip 5mm\cup\langle X_{\alpha_1},
X_{\alpha_1+\alpha_2}, X_{-\alpha_2}\rangle\cup
\langle X_{\alpha_1}, X_{-\alpha_2},
X_{-\alpha_1-2\alpha_2}\rangle\\
&\hskip 5mm\cup\langle X_{\alpha_1}, X_{-\alpha_2}+
\sqrt{-1}X_{\alpha_1+\alpha_2},
X_{\alpha_1+2\alpha_2}\rangle\\
&\hskip 5mm\cup \langle X_{\alpha_1}, X_{-\alpha_2}
-\sqrt{-1}X_{\alpha_1+\alpha_2},
X_{\alpha_1+2\alpha_2}\rangle.
\end{aligned}
\end{equation}

\noindent Therefore, by \eqref{nonsubreg} and
\eqref{333so5} we have
$$
\cod_{{\mathfrak a}_{X_{\alpha_1}}}\mathcal
I({\mathfrak a}_{X_{\alpha_1}})= 3.
$$
\end{nothing*}

\begin{nothing*}\label{G2}
{\it Case ${\mathfrak a}={\rm Lie}\,{\bf G}_2$.}

It is easy to verify   (see also \cite[p.\,10]{GQT})
that
\begin{equation*}
e=X_{\alpha_2}+X_{3\alpha_1+\alpha_2},\;
f=-2X_{-\alpha_2}-2X_{-3\alpha_1-\alpha_2},\; h=
2H_{\alpha_2}+2H_{3\alpha_1+ \alpha_2}\end{equation*}
is an ${\mathfrak s}{\mathfrak l}_2$-triple. Since
$\alpha_1(h)=0$, $\alpha_2(h)=2$, the classification
of nilpotent orbits in ${\mathfrak a}$,
see,\,e.g.,\,\cite[8.4]{CM}, implies that
\begin{equation}\label{xe}
x=e=X_{\alpha_2}+X_{3\alpha_1+\alpha_2}
\end{equation} is a subregular nilpotent element of
${\mathfrak a}$ and
\begin{equation}\label{subrG2}
\dim {\mathfrak a}_x=4.\end{equation}

Since $[X_\alpha, X_\beta]=0$ if $\alpha+\beta\notin
\Phi$, we have $X_{\alpha_1+\alpha_2},
X_{3\alpha_1+2\alpha_2}, X_{2\alpha_1+\alpha_2}\in
{\mathfrak a}_x$. Hence it follows from \eqref{xe},
\eqref{subrG2} that
\begin{equation}\label{centralizerg2}
{\mathfrak a}_x=\langle
X_{\alpha_2}+X_{3\alpha_1+\alpha_2},
X_{\alpha_1+\alpha_2}, X_{3\alpha_1+2\alpha_2},
X_{2\alpha_1+\alpha_2}\rangle.
\end{equation}
From \eqref{xe}, \eqref{centralizerg2} we deduce
that, for every $y\in {\mathfrak a}_x$, we have $[y,
{\mathfrak a}_x]\in \langle
X_{3\alpha_1+2\alpha_2}\rangle$; whence
\begin{equation*}
\begin{aligned}
\dim ({\mathfrak a}_x)_{y}=\dim {\mathfrak a}_{x}\!\!-\dim[y, {\mathfrak a}_x
] \geqslant 4-1=3.
\end{aligned}
\end{equation*}
This proves that $\mathcal I\bigl({\mathfrak
a}_{x}\bigr)={\mathfrak a}_x$, i.e.,
\begin{equation*}
\cod_{{\mathfrak a}_{x}}\mathcal I({\mathfrak
a}_{x})= 0.
\end{equation*}
\end{nothing*}

\begin{nothing*} As for the algebras ${\mathfrak a}$ of rank 3,
it will be sufficient for our purposes
 to consider
 only those ${\mathfrak a}$
 whose simple
ideals are of type ${\sf A}$, and to prove that for
such ${\mathfrak a}$ the inequality $\cod_{{\mathfrak
a}_x}\mathcal I({\mathfrak a}_x)\geqslant 1$ always
holds. We deduce this statement from the following
general lemma.
\end{nothing*}

\begin{lemma}\label{AAA} Let $x$ be a nilpotent
element of $\;{\mathfrak m}:={\mathfrak s}{\mathfrak
l}_{n_1}\oplus\ldots\oplus{\mathfrak s}{\mathfrak
l}_{n_q}$. Then there is semisimple element
$h\in{\mathfrak m}$ such that $x$ is a regular
element of $\;{\mathfrak m}_h$.
\end{lemma}

\vskip 2mm

\noindent {\it Proof.} Clearly, it suffices to prove
this statement for $q=1$. Therefore, we now assume
that ${\mathfrak m}={\mathfrak s}{\mathfrak l}_n$. If
$x\in {\mathfrak m}^{\rm reg}$, then $h=0$. Now let
$x\in {\mathfrak m}^{\rm irr}$. Then by the Jordan
normal form theory we may (and shall) assume that
\begin{equation}\label{xxx}
x={\rm diag}(J_{d_1},\ldots,J_{d_s}),\hskip 3mm
J_i=\begin{bmatrix}0&1&0&\ldots&0\\
0&0&1&\ldots&0\\[-3pt]
\hdotsfor[1]{5} \\
0&0&0&\ldots&1\\
0&0&0&\ldots&0\end{bmatrix}\in {\rm
Mat}_{i\times i},
\hskip 2mm 
s\geqslant 2.
\end{equation}

Since $s\geqslant 2$, there are $a_1,\ldots, a_s\in
\Bbb Z$ such that
\begin{gather}
a_1>\cdots> a_s,\label{aaa1}\\
d_1a_1+\cdots+d_sa_s=0.\label{aaa2}
\end{gather}
By \eqref{aaa2} the semisimple matrix ${\rm
diag}(a_1I_{d_1},\ldots,a_sI_{d_s})$ lies in
${\mathfrak m}$. We claim that one can take
\begin{equation*}
h:={\rm diag}(a_1I_{d_1},\ldots,a_sI_{d_s}).
\end{equation*}
Indeed, from \eqref{aaa1} we deduce that
\begin{equation}\label{centralizer}
[{\mathfrak m}_h,{\mathfrak m}_h]=\{{\rm
diag}(A_1,\ldots, A_s)\mid A_i\in {\mathfrak
s}{\mathfrak l}_{d_i}\hskip 2mm \mbox{for every
$i$}\},
\end{equation}
and the claim readily follows from \eqref{xxx} and
\eqref{centralizer}.\quad $\square$

\vskip 2mm

\begin{corollary}\label{>0} Maintain the notation of
Lemma {\rm \ref{AAA}}. Then
$$
\cod_{{\mathfrak m}_x}\mathcal I({\mathfrak
m}_x)\geqslant 1.
$$
\end{corollary}

\vskip 2mm

\noindent {\it Proof.} Since the element $h$ from
Lemma {\rm \ref{AAA}} is semisimple, the algebra
${\mathfrak m}_h$ is reductive and its rank is equal
to that of ${\mathfrak m}$. As $x$ is a regular
element of ${\mathfrak m}_h$, this implies the
equality $\dim ({\mathfrak m}_h)_x={\rm
rk}\,{\mathfrak m}$. But $({\mathfrak
m}_h)_x=({\mathfrak m}_x)_h$. Therefore,
$h\notin\mathcal I({\mathfrak m}_x)$; whence the
claim. \quad $\square$

\section{\bf Proofs of Theorems
\ref{ircod} and \ref{rational1}}

\begin{nothing*} {\it Proof of
Theorem {\rm \ref{ircod}}.} Let  $I$ be a subset of
$\Delta$ such that $1\leqslant |I|\leqslant 3$ and
let $x$ be a nilpotent element of ${\mathfrak
s}(I)^{\rm irr}$. If $|I|=1$, then according to
Subsection \ref{I1} we have $c(I,x)$=4. Now
con\-sider the cases $|I|=2$ and $3$.

(a) Let ${\mathfrak g}$ be of type ${\sf A}_r$, ${\sf
D}_r$, ${\sf E}_6$, ${\sf E}_7$, or ${\sf E}_8$. Then
every simple ideal of ${\mathfrak s}(I)$ is of type
${\sf A}$. Therefore, if $|I|=2$, then ${\mathfrak
s}(I)$ is isomorphic to ${\mathfrak s}{\mathfrak
l}_2\oplus{\mathfrak s}{\mathfrak l}_2$ or
${\mathfrak s}{\mathfrak l}_3$; whence by \eqref{x=0}
and Subsections \ref{sl2sl2}, \ref{sl3} we have $c(I,
x)=5$. If $|I|=3$, then ${\mathfrak s}(I)$ is
isomorphic to ${\mathfrak s}{\mathfrak
l}_2\oplus{\mathfrak s}{\mathfrak
l}_2\oplus{\mathfrak s}{\mathfrak l}_2$, ${\mathfrak
s}{\mathfrak l}_2\oplus{\mathfrak s}{\mathfrak l}_3$,
or ${\mathfrak s}{\mathfrak l}_4$;
whence by
Corollary \ref{>0} we have $c(I, x)\geqslant 4$. As
is explained in Subsection \ref{explain}, this
 information
implies  that $\cod_{\mathcal C}\irC=4$.

(b) Let ${\mathfrak g}$ be of type ${\sf B}_r$, ${\sf
C}_r$, or ${\sf F}_4$. If $|I|=2$, then ${\mathfrak
s}(I)$ is isomorphic to ${\mathfrak s}{\mathfrak
o}_5$, ${\mathfrak s}{\mathfrak l}_2\oplus{\mathfrak
s}{\mathfrak l}_2$, or ${\mathfrak s}{\mathfrak
l}_3$, and there is $I$ such that ${\mathfrak
s}{\mathfrak l}(I)$ is isomorphic to ${\mathfrak
s}{\mathfrak o}_5$. By \eqref{x=0} and Subsections
\ref{so51}, \ref{so52}, if ${\mathfrak s}(I)$ is
isomorphic to ${\mathfrak s}{\mathfrak o}_5$, then
$c(I, x)\geqslant 3$ and there is $x$ such that $c(I,
x)=3$. On the other hand, as we have seen in (a), if
${\mathfrak s}(I)$ is isomorphic to ${\mathfrak
s}{\mathfrak l}_2\oplus{\mathfrak s}{\mathfrak l}_2$
or ${\mathfrak s}{\mathfrak l}_3$, then $c(I, x)=5$.
According to Subsection \ref{explain}, this
 information
implies  that $\cod_{\mathcal C}\irC=3$.

(c) Let ${\mathfrak g}$ be of type ${\sf G}_2$. Then
$I=\Delta$ and ${\mathfrak s}(I)={\mathfrak g}$.
According to Subsection \ref{G2}, there is $x$ such
that $c(I, x)=2$. As is explained in Subsection
\ref{explain}, this implies that $\cod_{\mathcal
C}\irC=2$. \quad $\square$
\end{nothing*}

\begin{nothing*}\label{proof112} {\it Proof of
Theorem {\rm \ref{rational1}}.} Let $N$ be the
normalizer of ${\mathfrak t}$ in $G$. We endow
${\mathfrak t}\oplus{\mathfrak t}$ with the natural
$N$-module structure. Let $G\times^{N} \!({\mathfrak
t}\oplus{\mathfrak t})$ be the algebraic homogeneous
vector $G$-bundle over $G/N$ with fiber ${\mathfrak
t}\oplus{\mathfrak t}$, see \cite[\S2]{Se},
\cite[4.8]{PV2}, \cite[2.17]{LPR}. Denote by $g\ast
t$ the image of point $(g, t)\in G\times ({\mathfrak
t}\oplus{\mathfrak t})$ under the natural projection
$G\times ({\mathfrak t}\oplus{\mathfrak
t})\rightarrow G\times^{N} \!({\mathfrak
t}\oplus{\mathfrak t})$.

Since ${\mathcal C}=\overline{G({\mathfrak
t}\times{\mathfrak t})}$, the natural $G$-equivariant
morphism
\begin{equation}\label{algvectbundle}
\varphi\colon  G\times^{N} \!\!({\mathfrak
t}\oplus{\mathfrak t})\rightarrow {\mathcal C},\quad
g\ast t\mapsto g(t),
\end{equation}
is dominant. We claim that $\varphi$ is a birational
isomorphism. As ${\rm char}\,k=0$, proving this claim
is equivalent to showing  that $\varphi^{-1} (x)$ is
a single point for points $x$ in general position in
${\mathcal C}$. To show that the latter property
holds, notice that there is a nonempty open
$G$-stable subset $U$ of ${\mathcal C}$ laying in
$G({\mathfrak t}^{\rm reg}\times{\mathfrak t}^{\rm
reg})$. Let $x\in U$. Since $\varphi$ is
$G$-equivariant, the fibers $\varphi^{-1}(x)$ and
$\varphi^{-1}(g(x))$ for $g\in G$ are isomorphic.
Hence it would be sufficient to prove that
$\varphi^{-1}(x)$ is a single point for $x\in
{\mathfrak t}^{\rm reg}\times{\mathfrak t}^{\rm
reg}$.

To do this, take a point $x=(x_1, x_2)\in {\mathfrak
t}^{\rm reg}\times {\mathfrak t}^{\rm reg}$. Since
$x_i\in {\mathfrak t}^{\rm reg}$, we have
\begin{equation}\label{zzz}
{\mathfrak g}_{x_i}={\mathfrak t}. \end{equation} Let
$g\ast t\in \varphi^{-1}(x)$ where $t=(t_1, t_2)\in
{\mathfrak t}\oplus{\mathfrak t}$. By
\eqref{algvectbundle} we have $g(t_i)=x_i$. This and
\eqref{zzz} yield
\begin{equation}\label{=}
g({\mathfrak g}_{t_i})={\mathfrak
g}_{g(t_i)}={\mathfrak g}_{x_i}={\mathfrak t}.
\end{equation}
It follows from \eqref{=} that $\dim {\mathfrak
g}_{t_i}=\dim\,{\mathfrak t}$. Since ${\mathfrak
t}\subseteq {\mathfrak g}_{t_i}$, this yields
${\mathfrak g}_{t_i}={\mathfrak t}$. By virtue of
\eqref{=} the latter equality implies that $g\in N$.
From this and the definition of $g\ast t$ we now
deduce that $g\ast t=e\ast g(t)=e\ast x$. Thus,
$\varphi^{-1}(x)=e\ast x$. This proves the claim.

Since algebraic homogeneous vector bundles are
locally trivial in Zariski topo\-lo\-gy, see
\cite{Se}, the varieties $G\times^{N} \!\!({\mathfrak
t}\oplus{\mathfrak t})$ and $G/N\times ({\mathfrak
t}\oplus{\mathfrak t})$ are birationally isomorphic.
Hence the varieties ${\mathcal C}$ and $G/N\times
({\mathfrak t}\oplus{\mathfrak t})$ are birationally
isomorphic
 as well. Given this, the proof comes to a close
 because
 the variety $G/N$ of
 maximal tori of $G$ is rational,
see \cite{C}, \cite[6.1]{G}, \cite[7.9]{BS}. \quad
$\square$
\end{nothing*}

\end{document}